\documentclass[final,2p,times,12pt]{elsarticle}
\usepackage{titlesec}

\titleformat{\paragraph}
{\normalfont\normalsize\bfseries}{\theparagraph}{1em}{}
\titlespacing*{\paragraph}
{0pt}{3.25ex plus 1ex minus .2ex}{1.5ex plus .2ex}
\usepackage[centertags]{amsmath}
\usepackage{amsfonts,amssymb,amsthm,newlfont,graphicx,subfigure,algorithmic,algorithm,fixltx2e,booktabs,bm}
\usepackage[numbers]{natbib}
\usepackage{epstopdf}
\usepackage{empheq}
\usepackage{algorithm}
\usepackage{pifont}
\usepackage{multirow}
\usepackage{natbib}
\usepackage{enumitem}
\usepackage{enumerate}
\numberwithin{algorithm}{section}

\usepackage[colorlinks=true]{hyperref}
\biboptions{comma,square,sort&compress}

\numberwithin{equation}{section}
\renewcommand{\theequation}{\thesection.\arabic{equation}}

\def\simgt{\,\hbox{\lower0.6ex\hbox{$>$}\llap{\raise0.3ex\hbox{$\sim$}}}\,}
\def\simlt{\,\hbox{\lower0.6ex\hbox{$<$}\llap{\raise0.3ex\hbox{$\sim$}}}\,}
\def\simgteq{\,\hbox{\lower0.6ex\hbox{$\ge$}\llap{\raise0.6ex\hbox{$\sim$}}}\,}
\def\simlteq{\,\hbox{\lower0.6ex\hbox{$\le$}\llap{\raise0.6ex\hbox{$\sim$}}}\,}
\makeatletter
\def\user@resume{resume}
\def\user@intermezzo{intermezzo}
\newcounter{previousequation}
\newcounter{lastsubequation}
\newcounter{savedparentequation}
\usepackage{lipsum}
\makeatletter
\def\ps@pprintTitle{%
 \let\@oddhead\@empty
 \let\@evenhead\@empty
 \def\@oddfoot{}%
 \let\@evenfoot\@oddfoot}
\makeatother

\def\simgt{\,\hbox{\lower0.6ex\hbox{$>$}\llap{\raise0.4ex\hbox{$\sim$}}}\,}
\def\simlt{\,\hbox{\lower0.6ex\hbox{$<$}\llap{\raise0.4ex\hbox{$\sim$}}}\,}
\def\simgteq{\,\hbox{\lower0.6ex\hbox{$\ge$}\llap{\raise0.6ex\hbox{$\sim$}}}\,}
\def\simlteq{\,\hbox{\lower0.6ex\hbox{$\le$}\llap{\raise0.6ex\hbox{$\sim$}}}\,}
\begin{document}

\begin{frontmatter}
%-----------------------------------------------------------------
\title{Generalized Baer and Generalized Quasi-Baer Rings of Skew Generalized Power Series}

\author[a]{R. M. Salem}
\ead{rsalem_02@ hotmail.com}

\author[a]{R. E. Abdel-Khalek}
\ead{ramy_ama@yahoo.com}

\author[b]{M. M. Hamam}
\ead{mostafahamam154@gmail.com}

\address[a]{Department of Mathematics, Faculty of Science, Al-Azhar Univ., Nasr City 11884, Cairo, Egypt.}
\address[b]{Department of Mathematics, Faculty of Science, Assiut Univ., Assiut 71515, Egypt.}
%%%%%%%%%%%%%%%%%%%%%%%%%%%%%%%%%%%%%%%%%%%%%%%%%%%%%%%%%%%%%%%%%%%%%%%%%%%%%%%%%%%%%%%
\begin{abstract}

\par Let $R$ be a ring with identity, $(S,\leq)$ an ordered monoid, $\omega:S \to End(R)$ a monoid homomorphism, and $A= R\left[\left[S,\omega \right]\right]$ the ring of skew generalized power series. The concepts of generalized Baer and generalized quasi-Baer rings are generalization of Baer and quasi-Baer rings, respectively. A ring  $R$  is called generalized right Baer (generalized right quasi-Baer) if for any non-empty subset $S$ (right ideal $I$) of $R$, the right annihilator of  $S^n \hspace{0.1cm}(I^n)$ is generated by an idempotent for some positive integer $n$. Left cases may be defined analogously. A ring $R$ is called generalized Baer (generalized quasi-Baer) if it is both generalized right and left Baer (generalized right and left quasi-Baer) ring. In this paper, we examine the behavior of a skew generalized power series ring over a generalized right Baer (generalized right quasi-Baer) ring and prove that, under specific conditions, the ring $A$ is generalized right Baer (generalized right quasi-Baer) if and only if  $R$ is a generalized right Baer (generalized right quasi-Baer) ring.
\end{abstract}
\end{frontmatter}
%%%%%%%%%%%%%%%%%%%%%%%%%%%%%%%%%%%%%%%%%%%%%%%%%%%%%%%%%%%%%%%%%%%%%%%%
%%%%%%%%                    Introduction                         %%%%%%%
%%%%%%%%%%%%%%%%%%%%%%%%%%%%%%%%%%%%%%%%%%%%%%%%%%%%%%%%%%%%%%%%%%%%%%%%
%\setcounter{equation}{0}
\section{Introduction}\label{sec1}

Throughout this article, $R$ denotes an associative ring with identity, and $r_R\left(S\right)=\left\{a\in R \hspace{0.1cm } | \hspace{0.1cm} sa=0, for \hspace{0.1cm} all \hspace{0.1cm} s\in S\right\}$ is the right annihilator of a nonempty subset $S$ in $R$. In \cite{4kaplansky1955rings}, Kaplansky introduced Baer rings as rings in which the right annihilator of every nonempty subset of $R$ is generated by an idempotent. Clark defined quasi-Baer rings in \cite{2clark1967twisted} as rings in which the right annihilator of every right ideal of $R$ is generated by an idempotent. Baer rings are clearly quasi-Baer rings. In a reduced ring $R$, $R$ is Baer if and only if $R$ is quasi-Baer.  The definitions of Baer and quasi-Baer rings are left-right symmetric by \cite[Theorem 3]{4kaplansky1955rings} and \cite[Lemma 1]{2clark1967twisted}.

According to Moussavi et al. \cite{16moussavi2005generalized}, a ring $R$ is called generalized right quasi-Baer if for any right ideal $I$ of  $R$, the right annihilator of  $I^n$ is generated by an idempotent for some positive integer $n$, depending on $I$. The class of generalized right quasi-Baer rings includes the right quasi-Baer rings and is closed under direct product and also under some kinds of upper triangular matrix rings. Example (4.4) in \cite{16moussavi2005generalized} is an example of a generalized right quasi-Baer ring which is not generalized left quasi-Baer, and hence the definition of generalized quasi-Baer ring is not left-right symmetric.

In \cite{17paykan2015generalization} K. Paykan and A. Moussavi defined a generalized right Baer rings as rings  in which the right annihilator of $S^n$ is generated by an idempotent for some positive integer $n$, where $S$ is a non-empty subset of $R$ and  $S^n$ is a set that contains elements  $a_1 a_2… a_n$ such that $a_i\in S$ for  $1\leq i \leq n$. A ring is called generalized Baer if it is both generalized right and left Baer ring. Baer rings are clearly generalized right (left) Baer. Also, the class of generalized right (left) Baer rings is obviously included in the classes of generalized right (left) quasi Baer rings. Example (2.2) in \cite{17paykan2015generalization} shows that there are various classes of generalized quasi-Baer ring which are not generalized Baer. Also, there are rich classes of generalized right Baer rings which are not Baer (see \cite[Example 2.3]{17paykan2015generalization}).

In \cite{19hamam2024semi} we examine the behavior of a skew generalized power series ring over a semi-Baer (semi-quasi Baer) rings. In this paper, we study the relation between the generalized Baer (generalized quasi Baer) rings and its skew generalized power series ring extensions and determine the conditions under which a ring of skew generalized power series $ R\left[\left[S,\omega \right]\right]$ is generalized Baer (generalized quasi Baer) whenever $R$ is generalized Baer (generalized quasi Baer) and vice versa.

 %\newpage
%sssssssssssssssssssssssssssssssssssssssssssssssssssssssssssssssssssssssss
%\section{Preliminaries}\label{sec2}
\section{Skew Generalized Power Series Rings}\label{sec2}

The construction of generalized power series rings was considered by Higman in \cite{3higman1952ordering}. Paulo Ribenboim studied extensively in a series of  papers (see \cite{11ribenboim1991rings,12ribenboim1992noetherian,13ribenboim1994rings,14ribenboim1995special,15ribenboim1997semisimple})  the rings of generalized power series. In \cite{9mazurek2008neumann} Mazurek and Ziembowski generalized this construction by introducing the concept of the skew generalized power series rings.\\

An ordered monoid is a pair $(S,\leq)$ consisting of a monoid $S$ and a compatible order relation $\leq$ such that if  $u\leq v$, then  $ut\leq vt$ and $tu \leq tv$ for each $t\in S$. $(S,\leq)$ is called a strictly ordered monoid if whenever  $u,v \in S$  such that  $u < v$ (i.e., $u\leq v$ and $u\neq v$), then $ut < vt$  and  $tu < tv$  for all  $t\in S$. Recall that an ordered set $(S,\leq)$ is called artinian if every strictly decreasing sequence of elements of $S$ is finite, and $(S,\leq)$ is called narrow if every subset of pairwise order-incomparable elements of  $S$ is finite. Thus $(S,\leq)$  is artinian and narrow if and only if every nonempty subset of  $S$  has at least one but only a finite number of minimal elements.\\

Let $R$ be a ring,  $(S,\leq)$  a strictly ordered monoid, $\omega:S \to End(R)$   a monoid homomorphism, where $\omega_s$ denote the image of  $s$ under $\omega$, for each $s\in S$, that is  $\omega_s=\omega(s)$, and $A$ the set of all maps $f : S\to R$  such that  $supp (f) =\left\{s\in S :f(s)\neq 0\right\}$ is artinian and narrow subset of  $S$. Under pointwise addition $A$  is an abelian subgroup of the additive group of all mappings $f : S\to R$. For every $s\in S$ and $f,g \in A$  the set $X_s\left(f,g\right)=\left\{\left(u,v)\in S\times S : uv=s,f\left(u\right) \neq 0,g\left(v\right) \neq 0 \right)   \right\}$ is finite by \cite[4.1]{12ribenboim1992noetherian}. Define the multiplication for each $f,g\in A$ by:\\
 $fg\left(s\right)=\sum_{\left(u,v \right)\in X_s\left(f,g\right)}^{} f\left(u\right) \omega_u  (g\left(v\right))$. (by convention, a sum over the empty set is 0). With pointwise addition and multiplication as defined above, $A$  becomes  a ring called the ring of skew generalized power series whose elements have coefficients in $R$ and exponents in $S$. For each $r \in R$ and $s \in S$ one can associate the maps $c_r,e_s \in A$ defined by :

\[ c_r(x) =
  \begin{cases}
    r       & \quad \text{if } x = 1_s\\
    0  & \quad \text{otherwise}
  \end{cases}
,\quad  e_s(x) =
  \begin{cases}
    1_R       & \quad \text{if } x = s\\
    0  & \quad \text{otherwise}
  \end{cases}\]
It is clear that  $r\to c_r$  is a ring embedding of  $R$  into  $A$  and  $s \to e_s$
is a monoid embedding of  $S$  into  the multiplicative monoid of  $A$
and $e_s c_r = c_{\omega_s (r)} e_s$.
Moreover, the identity element of $A$  is a map $e : S\to R$  defined by $ e\left(1_S\right) =\left(1_R\right)$  and  $e\left(s\right)=0$  for each $s \in S \backslash \left\{1_s \right\} $.

Let $R$ be a ring and $\sigma$ an endomorphism of $R$. The construction of the skew generalized power series rings generalizes many classical ring constructions such as the skew polynomial rings $R[x,\sigma]$ if $S=N \cup \left\{0 \right\}$ and $\leq$ is the trivial order, skew power series rings $R[[x,\sigma]]$  if $S=N \cup\left\{0 \right\}$ and $\leq$ is the natural linear order, skew Laurent polynomial rings $R[x,x^{-1};\sigma]$ if $S=Z$ and $\leq$ is the trivial order where $\sigma$  is an automorphism of $R$, skew Laurent power series rings $R[[x,x^{-1};\sigma]]$ if $S=Z$ and $\leq$ is the natural linear order where $\sigma$ is an automorphism of $R$. Moreover, the ring of polynomials $R[x]$, the ring of power series $R[[x]]$, the ring of Laurent polynomials $R[x,x^{-1}]$, and the ring of Laurent power series $R[[x,x^{-1}]]$ are special cases of the skew generalized power series rings, if we consider $\sigma$  to be the identity map of $R$.
%\newpage
%sssssssssssssssssssssssssssssssssssssssssssssssssssssssssssssssssssssssss
\section{Main Results}\label{sec2}
An ordered monoid $(S,\leq)$ is called positively ordered if  $1$  is the minimal element of  $S$.\\\\
\textbf{Definition 3.1} (\cite{1annin2002associated}).
An endomorphism $\sigma$ of a ring $R$ is called compatible if for all  $a,b\in R$,  $ab=0$  if and only if  $a \sigma(b)=0$.\\\\
\textbf{Definition 3.2} (\cite{5krempa1996some}).
 An endomorphism $\sigma$ of a ring $R$ is called rigid if for every   $a\in R$,  $a \sigma(a)=0$  if and only if  $a=0$.\\\\
Let $R$ be a ring, $(S,\leq)$ a strictly ordered monoid, and $\omega:S \to End(R)$ a monoid homomorphism. As in \cite{8marks2010unified}, a ring  $R$ is $S$-compatible ($S$-rigid) if $\omega_s$ is compatible (rigid) for every $s\in S$.\\\\
\textbf{Definition 3.3} (\cite{7marks2009new}).
 An ordered monoid $(S,\leq)$ is said to be quasitotally ordered (and  $\leq$  is called a quasitotal order on  S ) if  $\leq$  can be refined to an order $\preceq$ with respect to which  $S$  is a strictly totally ordered monoid.\\
%%%%%%%%%%%%%%%%%%%%%%%%%%%%%

Recall that a ring $R$ is said to be $(S,\omega)$-Armendariz if whenever $fg=0$ for $f,g\in R[[S,\omega]]$, then $f(s).\omega_s (g(t))=0$  for all  $s,t\in S$ (see \cite[Definition 2.1]{8marks2010unified}).\\\\
 %%%%%%%%%%%%%%%%%%%%%%%%%%%%%
 \textbf{Proposition 3.4} (\cite[Proposition 4.10]{8marks2010unified}).
  Let $R$ be a ring, $(S,\leq)$ a strictly ordered monoid, and  $\omega:S \to End(R)$ a monoid homomorphism. Assume that $R$ is $(S,\omega)$-Armendariz. If $f$ is an idempotent of  $ R\left[\left[S,\omega\right]\right]$, then $f (1)$ is an idempotent of  $R$ and  $f = c f (1)$.\\\\
 %%%%%%%%%%%%%%%%%%%%%%%%%%%%%
 \textbf{Proposition 3.5}.
 Let $R$ be an $(S,\omega)$-Armendariz ring, $(S,\leq)$ a quasitotally ordered monoid, and  $\omega:S \to End(R)$ a monoid homomorphism. Set  $ A= R\left[\left[S,\omega\right]\right]$   the ring of skew generalized power series.\\
      (1)  If  $A$  is a generalized right Baer ring, then  $R$  is  a generalized right Baer ring.\\
      (2)  If  $R$ is an $S$-compatible ring and $A$  is a generalized right quasi-Baer ring, then  $R$  is  a generalized right quasi-Baer ring.\\\\
\textbf{Proof}.
 (1) Let  $X$  be a non-empty subset of $R$. Then  $B=\{c_x:x\in X\}$  is a non-empty subset of $A$. Since  $A$  is a generalized right Baer, there exists  $f\in A$  such that  $r_A (B^n )=fA$  with $f^2=f$. Proposition 3.4 implies that $f(1)$ is an idempotent element of  $R$. We want to prove that  $r_R (X^n )=f(1)R$.
Since $f\in r_A (B^n )$,  we have $(c_{x_1 } c_{x_2}… c_{x_n}) f=0$    for all   $c_{x_1} c_{x_2}… c_{x_n}\in B^n$ and  $x_1,x_2,…,x_n\in X$.
Thus $0=(c_{x_1} c_{x_2}… c_{x_n}) f (1)=c_{x_1} (1) \omega_1 (c_{x_2}(1))… \omega_1 (c_{x_n}(1)) \omega_1 (f(1))= x_1  x_2… x_n  f(1)$ for all   $x_1  x_2… x_n\in X^n$. Hence  $f(1)\in r_R (X^n )$,  which implies that  $f(1)R\subseteq r_R (X^n )$.
On the other hand, if $a\in r_R (X^n )$, then  $(x_1  x_2… x_n)a=0$  for all  $x_i\in X$  with  $1\leq i\leq n$. Thus $(c_{x_1} c_{x_2}… c_{x_n}) c_a  (1)=c_{x_1}(1) \omega_1 (c_{x_2 } (1))… \omega_1 (c_{x_n} (1))\omega_1 ( c_a (1))=(x_1  x_2… x_n) a=0$. Which implies that $(c_{x_1} c_{x_2}… c_{x_n}) c_a=0$  for all  $c_{x_i }\in B$. Therefore, $c_a\in r_A (B^n )=fA$  and  $c_a=fg$  for some  $g\in A$.  Now,  $a=c_a (1)=(fg)(1)=f(1) \omega_1 (g(1))\in f(1)R$. That is  $r_R (X^n )\subseteq f(1)R$, which follows that  $r_R (X^n )=f(1)R$. Hence  $R$  is  a generalized right Baer ring.\\

 (2) Let  $I$  be a right ideal of  $R$. Then $ I\left[\left[S,\omega\right]\right]=\{f \in A  |  f(s)\in I \enspace for\,any \enspace  s \in S\}$  is a right ideal of $A$. Since  $A$  is a generalized right quasi-Baer, there exists  $f\in A$  such that  $r_A (I^n [[S,\omega]])=fA$  with $f^2=f$. Proposition 3.4 implies that $f(1)$ is an idempotent element of  $R$. We want to prove that  $r_R (I^n )=f(1)R$. Since $f\in r_A (I^n [[S,\omega]])$,  we have $(g_1 g_2… g_n) f=0$    for all   $g_1,g_2,… ,g_n\in I[[S,\omega]]$.
Since  $c_{i_k }\in I[[S,\omega]]$  for all  $i_k\in I$  with  $1\leq k \leq n$, we have  $(c_{i_1 } c_{i_2}… c_{i_n }) f=0$. Consequently, $((c_{i_1 } c_{i_2}… c_{i_n }) f)(1)=c_{i_1} (1) \omega_1 (c_{i_2} (1))… \omega_1 (c_{i_n} (1)) \omega_1 (f(1))=0$  which implies that  $i_1  i_2… i_n  f(1)=0$ for all   $i_1,i_2,…,i_n\in I$. Hence  $f(1)\in r_R (I^n )$,  which implies that  $f(1)R\subseteq r_R (I^n )$.
On the other hand, if   $a\in r_R (I^n )$, then  $(i_1 i_2…i_n)a=0$  for all  $i_1,i_2,…,i_n \in I$. Since $g_k (s_k )\in I$  for all  $g_k\in I[[S,\omega]]$ and $ s_k \in S$  with  $1\leq k \leq n$, we have $g_1 (s_1 ) g_2 (s_2 )… g_n (s_n )  a=0$. Since $R$ is $S$-compatible, we have $g_1 (s_1 ) \omega_{s_1} (g_2 (s_2 )) \omega_{s_1 s_2 } (g_3 (s_3 ))…\omega_{s_1 s_2… s_{n-1} } (g_n (s_n )) \omega_{s_1 s_2… s_n } (c_a (1))=0$. \\ Which implies that $(g_1 g_2… g_n c_a )(s)=$\\
   $\sum_{\left(s_1,s_2,… , s_n,1\right)\in X_s\left(g_1,g_2,… ,g_n,c_a\right)}^{} g_1 (s_1 ) \omega_{s_1} (g_2 (s_2 )) \omega_{s_1 s_2 } (g_3 (s_3 ))… \omega_{s_1 s_2… s_n } (c_a (1))=0$.\\
Thus  $c_a\in r_A (I^n [[S,\omega]])=fA$  and  $c_a=fg$  for some  $g\in A$.  Now, $a=c_a (1)=(fg)(1)=
f(1) \omega_1 (g(1))\in f(1)R$. That is  $r_R (I^n )\subseteq f(1)R$, which follows that $r_R (I^n )=f(1)R$. Hence  $R$  is  a generalized  right quasi-Baer ring.\\\\
%%%%%%%%%%%%%%%%%%%%%%%%%%%%%
\textbf{Proposition 3.6}.
Let $R$  be  an  $S$-compatible $(S,\omega)$-Armendariz ring, $(S,\leq)$ a quasitotally ordered monoid and  $\omega:S \to End(R)$  a monoid homomorphism. Set  $ A= R\left[\left[S,\omega\right]\right]$  the ring of skew generalized power series.\\
      (1)  If  $R$ is a generalized right Baer ring, then $A$  is a generalized right Baer ring.\\
      (2) If  $R$ is a generalized right quasi-Baer ring, then $A$  is a generalized right quasi-Baer ring.\\\\
\textbf{Proof}.
(1) Let $B$ be a non-empty subset of  $A$. Then  $U= \{ f(s): f \in B ,s \in S \}$ is a non-empty subset of $R$. Since $R$ is a generalized right Baer, there exists  $b\in R$   such that  $r_R (U^n )=bR$ with $b^2= b$  which implies that $c_b^2= c_b$. We want to prove that  $r_A (B^n )=c_b A$.
Since $b\in r_R (U^n )$, it follows that   $ f_1 (s_1 ) f_2 (s_2 )…f_n (s_n )  b=0$  for all  $f_i (s_i )\in U$  with  $1\leq i\leq n$. Thus  $f_1 (s_1 ) f_2 (s_2 )…f_n (s_n )  c_b (1)=0$.  Since $R$ is $S$-compatible, then $f_1 (s_1 ) \omega _{s_1} (f_2 (s_2 ))…\omega_{s_{n-1} } (f_n (s_n ))\omega_{s_n} (c_b (1))=0$.
Thus   $(f_1 f_2… f_n  c_b )(s)=$\\
$\sum_{\left(s_1,s_2,… , s_n,1\right)\in X_s\left(f_1,f_2,…,f_n,c_b\right)}^{} f_1 (s_1 ) \omega_{s_1} (f_2 (s_2 )) \omega_{s_1 s_2 } (f_3 (s_3 ))… \omega_{s_1 s_2… s_n } (c_b (1))=0$.\\ It follows that  $c_b\in r_A (B^n )$  which implies that   $c_b A\subseteq r_A (B^n )$.\\
Now, let  $f\in r_A (B^n )$.   Then   $f_1 f_2… f_n  f=0$   for all  $f_1 f_2… f_n\in B^n$. Since $R$ is an  $(S,\omega)$-Armendariz ring, we get  $f_1 (u_1)\omega_{u_1 } (f_2 (u_2))… \omega_{u_{n-1} } (f_n (u_n ))\omega_{u_n } (f(v))=0$  for all  $u_1,u_2,…,u_n,v\in S$. Moreover, Since $R$ is $S$-compatible, we get \\ $f_1 (u_1 )  f_2 (u_2 )… f_n (u_n )  f(v)=0$.  Thus $f(v)\in r_R (U^n )=bR$  for all  $v\in S$. Therefore, for all $v\in S$ there exists  $r\in R$  such that  $f(v)=br=(c_b c_r e_v)(v)$. Thus  $f=c_b c_r e_v $ , which implies that  $f\in c_b A$. That is  $r_A (B^n )\subseteq c_b A$, which follows that  $r_A (B^n )=c_b A$. Hence  $A$  is a generalized right Baer ring.\\

(2) Let $J$  be a right ideal of  $A$. For every $s\in S$, set  $J_s=\{ f(s) | f \in J ,s \in S\}$, and  $J^*=\cup_{(s\in S)}  J_s$. Let  $I$  be the right ideal generated by $J^*$.  Since $R$ is a generalized right quasi-Baer ring, there exists $b\in R$ such that $r_R (I^n )=bR$ with  $b^2=b$. Therefore,  $c_b$ is an idempotent element of $A$. We want to prove that  $r_A (J^n )=c_b A$.
Since  $b\in r_R (I^n )$, it follows that  $i_1 i_2 i_3… i_n  b=0$  for all  $i_j\in I$  with  $1\leq j\leq n$. Since $g_i (s_i )\in I$  for all  $g_i\in J$  and  $s_i\in S$, we have  $g_1 (s_1 ) g_2 (s_2 )…g_n (s_n )  b=0$. Thus  $g_1 (s_1 ) g_2 (s_2 )…g_n (s_n )  c_b (1)=0$. Since $R$ is  $S$-compatibe,  $g_1 (s_1 ) \omega_{s_1 } (g_2 (s_2 )) \omega_{s_1 s_2 } (g_3 (s_3 ))…\omega_{s_1 s_2… s_{n-1}} (g_n (s_n )) \omega_{s_1 s_2… s_n} (c_b (1))=0$.  Thus $(g_1 g_2… g_n  c_b )(s)=$ \\
$\sum_{\left(s_1,s_2,… , s_n,1\right)\in X_s\left(g_1,g_2,…,g_n,c_b\right)}^{} g_1 (s_1 ) \omega_{s_1} (g_2 (s_2 )) \omega_{s_1 s_2 } (g_3 (s_3 ))… \omega_{s_1 s_2… s_n } (c_b (1))=0$.
It follows that  $c_b \in r_A (J^n )$  which implies that   $c_b A\subseteq r_A (J^n )$.\\
Now, let  $g \in r_A (J^n )$.   Then   $g_1 g_2… g_n  g=0$   for all  $g_1 ,g_2,…,g_n\in J$. Since $R$ is an $(S,\omega)$-Armendariz ring, we get  $g_1 (u_1)\omega_{u_1} (g_2 (u_2))… \omega_{u_{n-1}} (g_n (u_n ))\omega_{u_n } (g(v))=0$  for all  $u_1,u_2,…,u_n,v\in S$. Moreover, Since $R$ is $S$-compatible, we get \\ $g_1 (u_1 )  g_2 (u_2 )… g_n (u_n )  g(v)=0$. Thus  $g(v)\in r_R (I^n )=bR$  for all  $v\in S$. Therefore, for all  $v\in S$  there exists  $r\in R$  such that  $g(v)=br=(c_b c_r e_v)(v)$. Thus  $g=c_b c_r e_v$  , which implies that  $g\in c_b A$. That is  $r_A (J^n )\subseteq c_b A$, which follows that  $r_A (J^n )=c_b A$. Hence  $A$  is a generalized right quasi-Baer ring.\\
%%%%%%%%%%%%%%%%%%%%%%%%%%%%%

By combining Proposition 3.5 and Proposition 3.6, we obtain the following Theorem.\\\\
%%%%%%%%%%%%%%%%%%%%%%%%%%%%%
\textbf{Theorem 3.7}.
Let $R$ be an  $S$-compatible $(S,\omega)$-Armendariz ring, $(S,\leq)$  a quasitotally ordered monoid and $\omega:S \to End(R)$ a monoid homomorphism. Set $ A= R\left[\left[S,\omega\right]\right]$   the ring of skew generalized power series. Then $A$   is a generalized right Baer (quasi-Baer) ring  if and only if  $R$  is a generalized right Baer (quasi-Baer) ring.\\
%%%%%%%%%%%%%%%%%%%%%%%%%%%%%

Liu Zhongkui called a ring  $R$ an $S$-Armendariz ring if whenever $f,g\in R[[S]]$ (the ring of generalized power series) satisfy  $fg=0$, then $f(u) g(v)=0$ for each  $u,v\in S$ (see \cite{6liu2004special}).\\\\
\textbf{Corollary 3.8}.
Let  $R$  be  an  $S$-Armendariz ring and $(S,\leq)$  a quasitotally ordered monoid. Set  $A=R[[S]]$  the ring of generalized power series. Then  $A$   is a generalized right Baer (quasi-Baer) ring  if and only if  $R$  is a generalized right Baer (quasi-Baer) ring.\\
%%%%%%%%%%%%%%%%%%%%%%%%%%%%%

From \cite{20kim2006power}, a ring $R$ is called a power-serieswise Armendariz ring if whenever power series $f(x)= \sum_{ i=0 }^{\infty} a_i x^i$ and $g(x)= \sum_{ j=0 }^{\infty} b_i x^j$  satisfy  $f(x)g(x)=0$  we have  $a_i b_j=0$ for every  $i$ and $j$.\\\\
\textbf{Corollary 3.9}.
Let  $R$  be  a power-serieswise Armendariz ring. Then $ R[[x]]$  is a generalized right quasi-Baer ring  if and only if  $R$  is a generalized right quasi-Baer ring.\\\\
%%%%%%%%%%%%%%%%%%%%%%%%%%%%%
\textbf{Corollary 3.10} (\cite[Theorem 3.20 and Theorem 3.21]{17paykan2015generalization}).
Let  $R$  be  a power-serieswise Armendariz ring. Then  $R[[x]]$  is a generalized right Baer ring  if and only if  $R$  is a generalized right Baer ring.\\

Rege and Chhawchharia in \cite{10rege1997armendariz} introduced the notion of an Armendariz ring. They defined a ring $R$ to be an Armendariz ring if whenever polynomials $f(x)= \sum_{ i=0 }^{m} a_i x^i$ , $g(x)= \sum_{ j=0 }^{n} b_j x^j$ $\in R[x]$ satisfy  $f(x)g(x)=0$, then $a_i b_j=0$ for every  $i$ and $j$. (The converse is always true.) The name ‘‘Armendariz ring’’ was chosen because Armendariz \cite[Lemma 1]{21armendariz1974note} had noted that a reduced ring satisfies this condition. Note that Power-serieswise Armendariz rings are Armendariz, however the converse need not be true by example (2.1) in \cite{20kim2006power}.\\\\
%%%%%%%%%%%%%%%%%%%%%%%%%%%%%
\textbf{Corollary 3.11} (\cite[Proposition 1 and Proposition 2]{18javadi2010polynomial}).
Let  $R$  be  an Armendariz ring. Then  $R[x]$  is a generalized right quasi-Baer  ring  if and only if  $R$  is  a generalized right quasi-Baer ring.\\\\
%%%%%%%%%%%%%%%%%%%%%%%%%%%%%
\textbf{Corollary 3.12} (\cite[Theorem 3.14 and Theorem 3.15]{17paykan2015generalization}).
Let  $R$  be  an Armendariz ring. Then  $R[x]$  is a generalized right Baer  ring  if and only if  $R$  is  a generalized right Baer ring.\\\\
%%%%%%%%%%%%%%%%%%%%%%%%%%%%%

%-------------------------
%\bibliographystyle{abbrv}
%\bibliographystyle{alpha}
%\bibliographystyle{acm}
%\bibliography{\jobname}
%\bibliographystyle{apalike}
%\bibliographystyle{ieeetr}
%\bibliographystyle{plain}
\bibliographystyle{siam}
\bibliography{ref}

%=================================================================================
\end{document}